\documentclass{cpaper}

\usepackage{amsfonts, amsthm}
\newcommand{\Z}{\mathbb{Z}}
\newcommand{\Q}{\mathbb{Q}}

\input amssym.def
\newsymbol\square 1003
\newsymbol\lozenge 1006
\newsymbol\surjarrow 1310
\newsymbol\injarrow 131A
\newsymbol\ltimes 226E
\newsymbol\rtimes 226F
\newsymbol\void 203F
\newsymbol\leftrightsquigarrow 1321

\def\proclaim #1. #2\par #3\par {\medbreak
\noindent{\bf#1.\enspace}{\sl#2}\par\medbreak
\noindent{\bf Proof.} #3 \par}
\def\proclaimb #1. #2\par {\medbreak
\noindent{\bf#1.\enspace}{\sl#2}\par\medbreak}
 0
 2
 0
\newtheorem{thm}[equation]{Theorem}
\newtheorem{lemma}[equation]{Lemma}

\newtheorem{cor}[equation]{Corollary}

\begin{document}
\title{Stable rational cohomology of automorphism groups of free groups and the integral cohomology of moduli spaces of graphs}
\author{Craig A. Jensen}
\date{October 15, 2001}

\smallskip

\address{Department of Mathematics, University of New Orleans\\
New Orleans, LA 70148, USA} 
\email{jensen@math.uno.edu}

\begin{abstract}
It is not known whether or not the stable rational cohomology groups \newline
$\tilde H^*(Aut(F_\infty);\Q)$ always vanish (see Hatcher in \cite{[H]} and
Hatcher and Vogtmann in \cite{[V]} where they pose the question
and show that it does vanish in the first 6 dimensions.)
We show that either the rational cohomology does not vanish in
certain dimensions, or the integral cohomology of a moduli
space of pointed graphs does not stabilize in certain other dimensions.
Similar results are stated for groups of outer automorphisms.
This yields that $H^5(\hat Q_m; \mathbb{Z})$, 
$H^6(\hat Q_m; \mathbb{Z})$, and  $H^5(Q_m; \mathbb{Z})$
never stabilize as $m \to \infty$, where the moduli spaces
$\hat Q_m$ and $Q_m$ 
are the quotients of the spines $\hat X_m$ and $X_m$  of 
``outer space'' and ``auter space'', respectively, introduced 
in \cite{[C-V]} by Culler and Vogtmann and \cite{[H-V]} by Hatcher 
and Vogtmann.
\end{abstract}


\maketitlepage

\section{Introduction}

Let $F_n$ denote the free group on $n$ letters and let
$Aut(F_n)$ and $Out(F_n)$ denote the automorphism group
and outer automorphism group, respectively, of $F_n$.
In \cite{[H]} Hatcher shows that the integral cohomology
of the infinite symmetric group $\Sigma_\infty$ is a 
direct summand of the integral cohomology of $Aut(F_\infty)$.
He mentions that it is unknown whether or not the 
complementary summand is zero and in particular whether or
not $\tilde H^*(Aut(F_\infty);\Q)$ is always zero.  In
\cite{[H-V]}, Hatcher and Vogtmann again pose the question
of whether or not the stable rational cohomology groups of
$Aut(F_n)$ and $Out(F_n)$ all vanish, and show that it
does vanish in dimensions 1 through 6. A recent theorem of Madsen
and Tillman gives (after inverting the prime 2) 
a product decomposition for the plus
construction $B\Gamma^+$ of the classifying space for stable
mapping class groups; however, it is currently unknown to 
what extent this enables one to answer the question posed by
Hatcher and Vogtmann. 

Let $\hat X_m$ be the spine of outer space (see Culler and Vogtmann in 
\cite{[C-V]}) and let $\hat Q_m = \hat X_m / Out(F_m)$ be the
corresponding moduli space of graphs.  Similarly, let $X_m$ 
be the spine of auter space (see Hatcher and Vogtmann in 
\cite{[H-V]}) and let $Q_m = X_m / Aut(F_m)$ be the
corresponding moduli space of pointed graphs. 
In this paper, we show that

\begin{thm} \label{t46} 
Let $i \in \{0,1\}$. For all positive integers
$k$, either
$$H^{4k+i}(Out(F_\infty); \Q) \not =0$$
or
$$H^{4k+i+1}(\hat Q_m; \Z)
\hbox{ never stabilizes as } m \to \infty.$$
\end{thm}

\begin{thm} \label{t44} For all positive integers
$k$, either
$$H^{4k}(Aut(F_\infty); \Q) \not =0$$
or
$$H^{4k+1}(Q_m; \Z)
\hbox{ never stabilizes as } m \to \infty.$$
\end{thm}

From calculations in \cite{[V]} that
$$H^4(Aut(F_\infty);\Q) = H^4(Out(F_\infty);\Q) = 
H^5(Aut(F_\infty);\Q)=0,$$ the above two theorems 
immediately show that

\begin{cor} \label{t47} The cohomology groups $H^5(\hat Q_m; \Z)$
and $H^6(\hat Q_m; \Z)$ never
stabilize as $m \to \infty$.
\end{cor}

\begin{cor} \label{t45} The cohomology group $H^5(Q_m; \Z)$ never
stabilizes as $m \to \infty$.
\end{cor}

The two corollaries are true because as $m$ increases, 
torsion from increasingly higher primes is
introduced in $H^5(\hat Q_m; \Z)$, $H^6(\hat Q_m; \Z)$, and
$H^5(Q_m; \Z)$.  There are natural inclusions $Q_m \injarrow Q_{m+1}$, 
and it is known \cite{[H-V]} that the induced map
$H^i(Q_{m+1}; \Q) \to H^i(Q_m; \Q)$  
is an isomorphism for $m > 3i/2$.
It is therefore important
to keep in mind that the above two corollaries only hold with
respect to integral cohomology.

A quick note about our notation is appropriate here.
In general, groups without any additional
structure will be written using multiplicative notation
(e.g., $\Z/p \times \Z/p \cong (\Z/p)^2$) but modules like cohomology
groups will be written using additive notation
(e.g., $\Z/p \oplus \Z/p \cong 2(\Z/p)$.) 

In section 2 we review the basics about outer and auter space, and in
section 3 we prove Theorem \ref{t46}.  Symmetry groups of graphs 
with $2p-1$ holes are discussed in section 4, which enables us to 
prove Theorem \ref{t44} in section 5.

This paper is based on a dissertation 
(see \cite{[JD]}, \cite{[J1]}) written while the author 
was a student of Karen Vogtmann at Cornell, and the author would like
to thank Prof. Vogtmann for her help and advice.
The author would also like to thank Henry Glover for
his helpful comments on this paper.

\section{Basics about spectral sequences and $Aut(F_n)$}

Let $G$ be a group acting cellularly on a
finite dimensional CW-complex $X$ such that the stabilizer
$stab_G(\delta)$ of every cell $\delta$ is finite and such
that the quotient of $X$ by $G$ is finite.  Further suppose that
for every cell $\delta$ of $X$, the group $stab_G(\delta)$ fixes
$\delta$ pointwise.  Let $M$ be a $G$-module. 
Recall (see \cite{[B]}) that the equivariant
cohomology groups of the $G$-complex $X$
with coefficients in $M$ are defined by
$$H_G^*(X; M) = H^*(G; C^*(X;M))$$
and that if in addition $X$ is contractible
(which will usually, but not always, be the case in this
paper) then
$$H_G^*(X; M) = H^*(G;M).$$

In \cite{[B]} a spectral sequence
\begin{equation} \label{e1}
\tilde E_1^{r,s} = \prod_{[\delta] \in \Delta_n^r}
H^s(stab(\delta); M) \Rightarrow H_G^{r+s}(X; M)
\end{equation}
is defined, where $[\delta]$ ranges over the set $\Delta_n^r$ of orbits
of $r$-simplices $\delta$ in $X$.

If $M$ is $\Z/p$ or $\Z_{(p)}$ then
a nice property should be noted about the spectral sequence \ref{e1}.
This property will greatly reduce the calculations we need to
go through, and in general will make concrete
computations possible.
Since each
group $stab(\delta)$ is finite, a standard restriction-transfer argument
in group cohomology yields that $|stab(\delta)|$ annihilates
$H^s(stab(\delta); M)$ for all $s > 0$.
(For examples of these sorts of
arguments see \cite{[B]}.)
Since all primes not equal to
$p$ are divisible in $\Z/p$ or $\Z_{(p)}$,
this in turn shows that the $p$-part of
$|stab(\delta)|$ annihilates $H^s(stab(\delta); M)$ for
$s>0$.
In particular, if $p$ does not divide some $|stab(\delta)|$, then this
$[\delta]$ does not contribute anything to the spectral sequence
\ref{e1} except in the horizontal row $s=0$.  It follows that if our
coefficients are $\Z/p$ or $\Z_{(p)}$ then we are mainly just concerned
with the simplices $\delta$ which have ``$p$-symmetry''.

We now specialize to the cases where
$G$ is $Out(F_n)$ or $Aut(F_n)$ and $X$ is either the spine 
$\hat X_n$ of ``outer space'' or the spine $X_n$ of ``auter space.''
Hatcher and Vogtmann's definition of auter space closely follows
Culler and Vogtmann's (prior) definition of outer space, except that
the graphs arising have basepoints.
We review some basic properties and definitions of
auter space below, where we concentrate on auter space
because that is where most of the calculations in this 
paper will take place.  Most of these facts can be found in \cite{[C-V]},
\cite{[H-V]}, \cite{[S-V1]}, and \cite{[S-V2]}.  

Consider the automorphism group $Aut(F_n)$ of
a free group $F_n$ of rank $n$ (where $n$ will be 
$2p-1$ for most of our work.)
Let $(R_n,v_0)$ be the
$n$-leafed rose,
a wedge of $n$ circles. We say a basepointed graph $(G,x_0)$
is \emph{admissible} if it has no free edges, all vertices except the
basepoint
have valence at least three, and there is a basepoint-preserving continuous
map $\phi\colon R_n \to G$ which induces an isomorphism on $\pi_1$. The
triple $(\phi,G,x_0)$ is called a \emph{marked graph}.  Two marked graphs
$(\phi_i,G_i,x_i) \hbox{ for } i=0,1$  are \emph{equivalent}  if there is
a homeomorphism $\alpha \colon (G_0,x_0) \to (G_1,x_1)$ such that
$ (\alpha\circ\phi_0)_\# = (\phi_1)_\# : \pi_1(R_n,v_0) \to \pi_1(G_1,x_1)$.
Define a partial order on the set of all equivalence classes of
marked graphs by setting $(\phi_0,G_0,x_0) \leq
(\phi_1,G_1,x_1)$ if $G_1$ contains a \emph{forest}
(a disjoint union of trees in $G_1$ which contains all of the vertices
of $G_1$) such that collapsing each tree in
the forest to a point yields $G_0$, where the collapse is compatible with
the maps $\phi_0$ and $\phi_1$.

From \cite{[H]} and \cite{[H-V]} we have that $Aut(F_n)$ acts
with finite stabilizers on
a contractible space $X_n$.
The space $X_n$ is the geometric realization of the poset of
marked graphs that we defined above.
Let $Q_n$ be the quotient of $X_n$ by
$Aut(F_n)$.
Note that the CW-complex $Q_n$ is not necessarily a simplicial
complex.
Since $Aut(F_n)$ has a torsion free subgroup of finite index
\cite{[H]} and it acts on the contractible,
finite dimensional space
$X_n$ with finite stabilizers and finite quotient,
$Aut(F_n)$ has finite vcd.

Let $p$ be an odd prime number, and let $\Z_{(p)}$ be the localization of
$\Z$ at the prime ideal $(p)$.  Then we
can apply the spectral sequence \ref{e1} to get
\begin{equation} \label{er1}
\tilde E_1^{r,s} = \prod_{[\delta] \in \Delta_n^r}
H^s(stab(\delta); \Z_{(p)}) \Rightarrow H^{r+s}(Aut(F_n); \Z_{(p)})
\end{equation}
where $[\delta]$ ranges over the set $\Delta_n^r$ of orbits
of $r$-simplices $\delta$ in $X_n$.

The spectral sequence \ref{er1}
requires as input the stabilizers $stab_{Aut(F_n)}(\delta)$ of
simplices $\delta$ in $X_n$.  Smillie and Vogtmann \cite{[S-V1]}
examined the structure of these stabilizers in
detail, and we list their results here.
Consider a given $r$-simplex
$$(\phi_r,G_r,x_r) > \cdots  > (\phi_1,G_1,x_1)
> (\phi_0,G_0,x_0)$$
with corresponding forest collapses
$$(H_r \subseteq G_r), \ldots, (H_2 \subseteq G_2), (H_1 \subseteq G_1).$$
For each $i \in {0,1, \ldots, r}$, let $F_i$ be the inverse image
under the map
$$G_r \to \cdots \to G_{i+1} \to G_i$$ of forest collapses,
of the
forest $H_i$.  That is, we
have
$$F_r \subseteq \cdots \subseteq F_2 \subseteq F_1 \subseteq G_r.$$
It is shown in \cite{[S-V1]}
that the stabilizer of the simplex under
consideration is isomorphic to the group
$Aut(G_r,F_1,\ldots,F_r,x_r)$ of basepointed automorphisms of the graph
$G_r$ that respect each of the forests $F_i$.  For example,
the stabilizer of a point $(\phi,G,x_0)$ in $X_n$ is isomorphic to
$Aut(G,x_0)$.

\section{Graphs without basepoints} \label{c15}

Theorem \ref{t46} is a direct consequence of the
stability theorems in \cite{[H]} and the spectral sequence
calculations in \cite{[G-M]}.

\begin{proof}[Proof of Theorem \ref{t46}:]
From \cite{[H]}, $$H^{4k+i}(Out(F_\infty); \Q) = H^{4k+i}(Aut(F_\infty); \Q)$$
and
if $m \geq 4k^2+10k+1+i^2/4+2ik+5i/2,$ then the standard map
$$H^{4k+i}(Aut(F_{m}); \Z) \to H^{4k+i}(Out(F_m); \Z)$$
is an isomorphism.
Observe that $H^{4k+i}(Out(F_{4k^2+10k+1+i^2/4+2ik+5i/2}); \Z) =
H^{4k+i}(Out(F_\infty); \Z)$ is a finitely generated
abelian group.  If it contains a torsion free summand isomorphic to
$\Z$, then we are done and $H^{4k+i}(Out(F_\infty); \Q) \not =0$.
Otherwise, choose a prime $q$ such that $q + 1 \geq 4k^2+10k+1+i^2/4+2ik+5i/2$
and so that for all primes $p \geq q$ there is no $p$-torsion
in $H^{4k+i}(Out(F_{\infty}); \Z).$
We will show that $H^{4k+i+1}(\hat Q_{p+1}; \Z)$ has $p$-torsion for
infinitely many primes $p$, which will prove the theorem.

Let $p \geq \hbox{ max}\{q, 25\}$ 
with $p \equiv 3 \hbox{ (mod } 4 \hbox{)}$.  (Note that there
are infinitely many possibilities for $p$, as there are 
infinitely many primes that are greater than a given 
number and congruent to $3$ modulo $4$.)
Because $H^{4k+i}(Out(F_{p+1}); \Z)$ has no
$p$-torsion, there is also no $p$-torsion in
$H^{4k+i}(Out(F_{p+1}); \Z_{(p)})$.
From the calculation of Glover and Mislin
in \cite{[G-M]} of the $E_2$-page of the
equivariant spectral sequence used to
calculate $H^{*}(Out(F_{p+1}); \Z_{(p)})$,
we know that this $E_2$-page,
in the rows $0 \leq s < 2(p-1)$,
is given by
\bigskip
$$\matrix{
\hfill E_2^{r,s} &=& &\left\{\matrix{
H^r(\hat Q_{p+1}; \Z_{(p)})  \hfill &s = 0 \hfill \cr
\Z/p \hfill &r=0 \hbox{ and } s = 4k > 0, k \in \Z^+ \hfill \cr
(n_p) \Z/p \hfill &r=1 \hbox{ and } s = 4k > 0, k \in \Z^+ \hfill \cr
0 \hfill &\hbox{otherwise} \hfill \cr} \right. \hfill \cr
}$$
\bigskip
where $n_p = (p-1)/12 - \epsilon_p$
and $\epsilon_p \in \{0,1\}$. Since $p \geq 25$, 
note that $n_p \geq 1$.

Hence a class
$\hat \alpha \in E_2^{i,4k}$ in the $E_2$-page survives at least
until the $E_{4k+1}$-page.
The class $\hat \alpha \in E_2^{i,4k}$ cannot survive to the
$E_\infty$ page, however, because there is no $p$-torsion in
the finite (since $H^{4k+i}(Out(F_{p+1}); \Q) = 0$)
additive group $H^{4k+i}(Out(F_{p+1}); \Z_{(p)})$.

It follows that there is $p$-torsion in
$$E_{4k+1}^{4k+i+1,0} = H^{4k+i+1}(\hat Q_{p+1}; \Z_{(p)}).$$
Thus $H^{4k+i+1}(\hat Q_{p+1}; \Z)$ has $p$-torsion. \end{proof}

\newpage

\section{Symmetry groups of graphs} \label{c13}

We will use spectral sequence \ref{er1} to compute a portion of 
the cohomology of $Aut(F_n)$.
Since our coefficient ring is $\Z_{(p)}$,
we have already remarked that
for the terms in the spectral sequence above the
horizontal axis, we are concerned only with simplices
whose stabilizers are divisible by $p$.
In addition, the stabilizer of a simplex consists of graph automorphisms
that respect the forest collapses in the simplex.  We will find
which simplices arise in the case $n=2p-1$.
In other words, we want to calculate which
graphs $G$ with a $\Z/p$
action on them have $\pi_1(G) \cong F_n$.
Recall that a
$\Z/p$-graph
$G$ is \emph{reduced} if it contains no $\Z/p$-invariant
subforests.

We now examine the cohomology of the quotient $Q_n$ of the
spine $X_n$ of auter space.  
There are natural inclusions $Q_m \injarrow Q_{m+1}$, and it is
known \cite{[H-V]} that the induced map
$H^i(Q_{m+1}; \Q) \to H^i(Q_m; \Q)$
is an isomorphism for $m > 3i/2$.
Our goal is to show that, in contrast,
$H^5(Q_m; \Z)$
never stabilizes as $m \to \infty$.  This is done by showing
that as $m$ increases, torsion from increasingly higher primes is
introduced in $H^5(Q_m; \Z)$.  
To this end, we do specific calculations in the spectral sequence
\ref{e1} applied to the action of $Aut(F_n)$ on $X_n$ for
$n=2p-1$.
The $E_2^{r,0}$-term of this spectral sequence is $H^r(Q_n; \Z_{(p)})$,
and the sequence converges to $H^r(Aut(F_n); \Z_{(p)})$.
Results from Hatcher and Vogtmann \cite{[H-V]} on the
cohomology of $Aut(F_n)$ are then used to obtain the result.

In this section, we do the ground work necessary to compute the 
$E_1$-page of the spectral sequence: we find all simplices of $X_n$ with 
$p$-symmetry and compute the cohomology of the stabilizers
of these simplices with coefficients in $\Z_{(p)}$.
In Section \ref{c14} we will compute the $E_2$-page of the
spectral sequence, and use this calculation to obtain the result.

Unless otherwise stated, $p \geq 5$ will be prime and $n=2p-1$.
The assumption that $p \geq 5$ is for convenience more than
any other reason, as the main results will only consider
arbitrarily large primes $p$ and so we should not
devote extra time to the (fairly easy to resolve) complications
introduced by considering the prime $p=3$. 
These complications arise from the fact that
the dihedral
group $D_6$ is the same as the symmetric group $S_3$,
so that we cannot
distinguish between dihedral and symmetric symmetry in that case.

\begin{figure}[here]
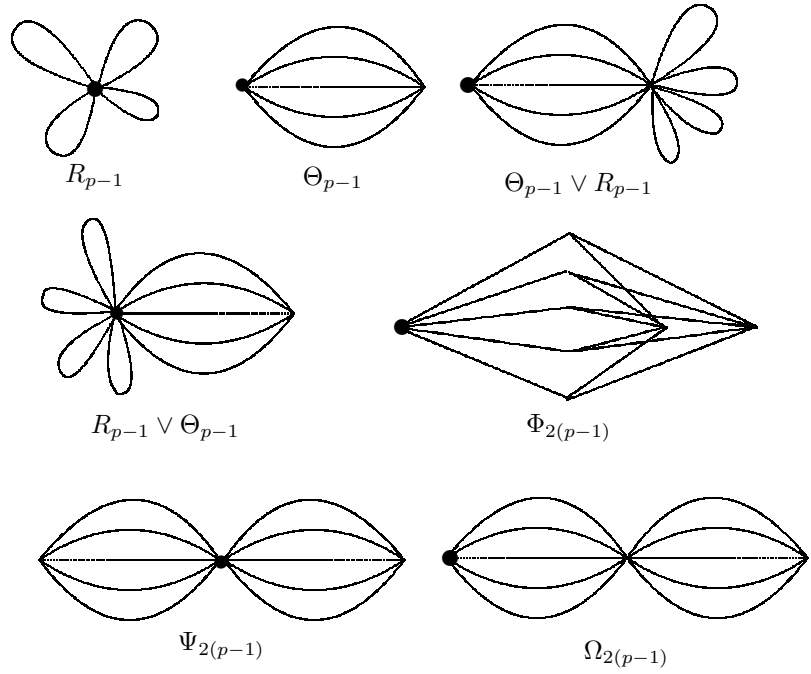

\caption{\label{fig1} Some graphs with $p$-symmetry}
\end{figure}

\begin{figure}[here]
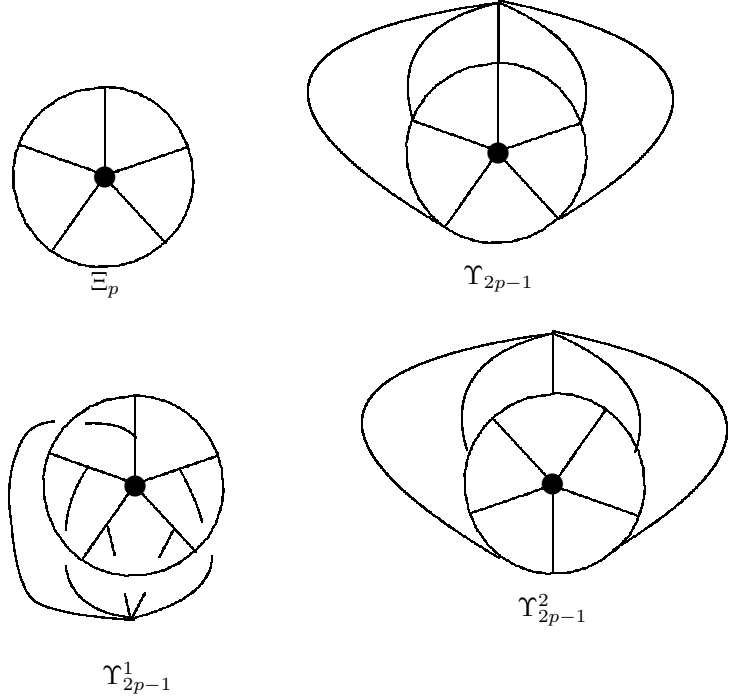

\caption{\label{fig2} Some graphs with $D_{2p}$-symmetry}
\end{figure}

We now define some graphs
that we will need for this section.
(Refer to Figures \ref{fig1} and \ref{fig2} for illustrations of most of these graphs.)
Let $\Theta_{p-1}$ be the graph with two vertices and $p$ edges, each of
which goes from one vertex to the other (see Figure \ref{fig1}.)
Say the ``leftmost vertex''
of $\Theta_{p-1}$ is the basepoint.  Hence when we write
$\Theta_{p-1} \vee R_{p-1}$ then we are stipulating that the rose
$R_{p-1}$ is attached to the non-basepointed vertex of $\Theta_{p-1}$,
while when we write $R_{p-1} \vee \Theta_{p-1}$ then we are saying
that the rose is attached to the basepoint of $\Theta_{p-1}$.
Let $\Phi_{2(p-1)}$ be a graph with $3p$ edges $a_1, \ldots, a_p,$
$b_1, \ldots, b_p,$ $c_1, \ldots, c_p,$ and $p+3$ vertices
$v_1, \ldots, v_p, x, y, z.$  The basepoint is $x$ and each of the
edges $a_i$ begin at $x$ and end at $v_i$.  The edges $b_i$ and $c_i$
begin at $y$ and $z$, respectively, and end at $v_i$.  Note that there
are obvious actions of $\Z/p$ on $\Theta_{p-1}$ and $\Phi_{2(p-1)}$,
given by rotation, and that these actions are unique up to conjugacy.
Let $\Psi_{2(p-1)}$ be the graph obtained from $\Phi_{2(p-1)}$ by
collapsing all of the edges $a_i$ to a point.  Let $\Omega_{2(p-1)}$
be the graph obtained from $\Phi_{2(p-1)}$ by collapsing either the
edges $b_i$ or the edges $c_i$ (the resulting graphs are isomorphic)
to a point.  Note that the only difference between $\Psi_{2(p-1)}$ and
$\Omega_{2(p-1)}$ is where the basepoint is located.

\input{fig1p.pic}

\newpage

Given a finite subgroup $G$  
of $Aut(F_n)$ for some integer $n$, we say that
a marked graph
$$\eta^1 : R_{r} \to \Gamma^1$$
is a {\em $G$-equivariant blowup in the fixed point space
$X_r^G$}
of a marked graph
$$\eta^2 : R_{r} \to \Gamma^2$$
if there is a $1$-simplex $\eta^1 > \eta^2$ in $X_r^G$.

\input{p4fig1.pic}

Let $\Upsilon_{2p-1}^1$ and $\Upsilon_{2p-1}^2$ be the two possible graphs
that can be obtained from $\Upsilon_{2p-1}$ by
equivariantly blowing
up 
the $p$ valence $4$ vertices into $2p$ valence $3$
vertices. 
That is, $\Upsilon_{2p-1}^1$ can be obtained by first taking a
$p$-gon and then attaching $p$ free edges to the $p$ vertices of
the $p$-gon.  Say each of these new edges $e_i$ begins at the vertex
$x_i$ and ends at the vertex $y_i$, and suppose that the vertices
$x_i$ are the ones that are attached to the $p$-gon.
Now form the $1$-skeleton of the
double cone or suspension over the $p$ vertices
$y_i$.  This gives the graph $\Upsilon_{2p-1}^1$.  The graph
$\Upsilon_{2p-1}^2$ can be thought of as follows:  First take a $p$-gon
and
cone off over the $p$ vertices of the $p$-gon.  Now also cone off
over the $p$ midpoints of the $p$ edges of the $p$-gon.
Note that there
is an obvious $\Z_p$-action on each of $\Theta_{p-1}$, $\Xi_p$,
$\Upsilon_{2p-1}$, $\Upsilon_{2p-1}^1$, and $\Upsilon_{2p-1}^2$.

Let $\Xi_p$ be the
$1$-skeleton of the cone over
a $p$-gon, so that $\Xi_p$ has $p+1$ vertices and $2p$ edges, 
one vertex has valence $p$ and the other $p$ vertices all have valence
$3$.  Let $\Upsilon_{2p-1}$ be the
$1$-skeleton of the 
suspension of a $p$-gon.  Hence
$\Upsilon_{2p-1}$ has $p+2$ vertices and $3p$ edges;  two of the vertices have
valence $p$ and the other $p$ have valence $4$.

\input{fil3p.pic}
\begin{figure}[here]
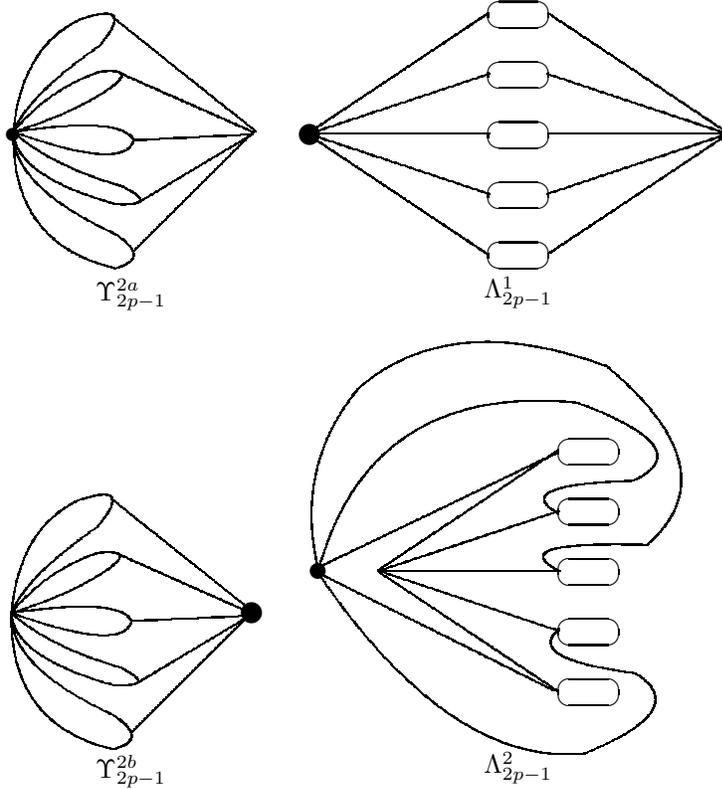

\caption{\label{fil3} Some graphs with (cohomologically) $\Sigma_p$-symmetry}
\end{figure}

Choose basepoints for the graphs $\Theta_{p-1}$, $\Upsilon_{2p-1}$,
$\Upsilon_{2p-1}^1$, and $\Upsilon_{2p-1}^2$ as
illustrated in Figure \ref{fig2}:  
Let the vertex on the ``leftmost'' side of $\Theta_{p-1}$ be
the basepoint.  Additionally, orient $\Upsilon_{2p-1}$,
$\Upsilon_{2p-1}^1$, and $\Upsilon_{2p-1}^2$
so that one of their
valence $p$ vertices is on the ``left'' and the other is on the
``right'' and choose the leftmost vertex to be the basepoint.
Writing $R_p \vee \Theta_{p-1}$ will mean that the two graphs are
wedged together at the basepoint of $\Theta_{p-1}$, while writing
$\Theta_{p-1} \vee R_p$ will mean that the non-basepointed vertex of
$\Theta_{p-1}$ is wedged to the vertex of $R_p$.  Let $\Upsilon_{2p-1}^{2a}$
be the graph obtained from $\Upsilon_{2p-1}^2$ by collapsing the
leftmost $p$ edges and let $\Upsilon_{2p-1}^{2b}$ be the one obtained
by collapsing the rightmost $p$ edges.
Refer to Figure \ref{fil3} for pictures of these graphs.
Figure \ref{fil3} also depicts two graphs
$\Lambda_{2p-1}^1$ and $\Lambda_{2p-1}^2$ which will be
used in the proof of Lemma \ref{t40}.

For the next lemma (Lemma \ref{t40})
only, we will consider the above graphs not to have
basepoints specified. 
The basepoint is just assumed to be located at some spot which is
invariant under the $\Z/p$ action. (This assumption is for
convenience rather than anything else,
so that we will not need to
introduce several separate subcases, each corresponding to
a different location for the basepoint.)

We will be looking at the standard equivariant spectral
sequence \ref{e1} applied to calculating the cohomology
groups $H^*(Aut(F_n); \Z_{(p)})$.  In particular, we
will be looking at the $E_1$ page of this
spectral sequence only in rows $0$ through $2p-3$, and often
just in rows $1$ though $2p-3$.  One interesting fact about these
rows is that they allow us to distinguish between simplices 
that have stabilizers whose cohomology is $\Z/p$, $D_{2p}$,
or $\Sigma_p$.
It is well
known that
$$\matrix{\hfill H^*(\Z/p; \Z_{(p)}) &=& \Z_{(p)}[x_2]/(px_2), \hfill \cr
\hfill H^*(D_{2p}; \Z_{(p)}) &=& \Z_{(p)}[x_4]/(px_4), \hfill \cr
\hfill \hbox{and } H^*(\Sigma_p; \Z_{(p)}) &=&
\Z_{(p)}[x_{2(p-1)}]/(px_{2(p-1)}), \hfill \cr}$$
where $x_2$, $x_4$ and $x_{2(p-1)}$ are generators of dimensions
$2$, $4$ and
$2(p-1)$, respectively.  Hence if a simplex of $X_n$ has
stabilizer isomorphic 
to $\Z/p$ or $D_{2p}$ then it will
contribute something to the $E_1$ page of the spectral
sequence in some of the rows $1$ through $2p-3$.  On the other hand,
if its stabilizer is isomorphic to $\Sigma_p$, then it
will not contribute anything to the $E_1$ page of the
spectral sequence in the given rows.

Define an $r$-simplex
in the $p$-singular locus of $X_n$ to have
{\em exactly $\Z/p$ symmetry} if it contributes exactly one copy of
$\Z/p$ to each of the entries $E_1^{r,2k}$, $0 < 2k < 2(p-1)$,
in the $E_1$ page of the spectral sequence.  Define an $r$-simplex
in the $p$-singular locus of $X_n$ to have {\em at most dihedral
symmetry}
if it contributes exactly one copy of
$\Z/p$ to each of the entries $E_1^{r,4k}$, $0 < 4k < 2(p-1)$,
in the $E_1$ page of the spectral sequence.

The next lemma examines which vertices in $X_n$ 
contribute to the spectral sequence in the given rows.
The proof of the
lemma actually explicitly enumerates which graphs have $p$-symmetry,
which will be very useful to us later.

\begin{lemma} \label{t40} Let $p \geq 5$ be prime and set $n=2p-1$.
Let the marked graph $(\xi,G,x_0)$ be a vertex in the
$p$-singular locus of $X_n$. 
Then the cohomology
$H^*(Aut(G,x_0); \Z_{(p)})$
of the stabilizer of this vertex is the
same as the cohomology with $\Z_{(p)}$ coefficients
of one of $D_{2p}$, $D_{2p} \times \Sigma_p$,
$\Sigma_p$, $\Sigma_p \times \Sigma_p$, 
$(\Sigma_p \times \Sigma_p) \rtimes \Z/2$
or $\Sigma_{2p}$.

\end{lemma}

\begin{proof} 
From \cite{[B-T]} and \cite{[M]},
we see that $p^2$ is an upper bound for the order of any $p$-subgroup
of $Aut(F_n)$.  Thus $p^2$ is an
upper bound for the order of a maximal $p$-subgroup $P$ of $Aut(G,x_0)$.
Since all possible choices for $P$ are abelian (i.e., there are only
three possibilities: $\Z/p$, $\Z/p^2$, and $\Z/p \times \Z/p$), we can
apply Swan's theorem (see \cite{[S]})
to see that
\begin{equation} \label{e8}
H^*(Aut(G,x_0); \Z_{(p)}) = H^*(P; \Z_{(p)})^{N_{Aut(G,x_0)}(P)}.
\end{equation}

We now look at each of the individual cases $P=\Z/p, \Z/p^2, \hbox{
and } \Z/p \times \Z/p$.

\medskip

\noindent CASE 1. We will first examine the case where
$P=\Z/p^2$.  In this case, we have that
$p^2$ edges $e_1, \ldots, e_{p^2}$ in $G$ are rotated around by $P$.
An examination of all possible ways that these edges could be
connected together, keeping in mind that $G$ is admissible, reveals
that this case is impossible.  For example, the first subcase is that
all of the $e_i$ begin and end at the same vertex.  This is not
possible because the fundamental group of $R_{p^2}$ is too large for
it to be a subgraph of $G$. For the next subcase, suppose each edge
goes from some vertex $y_1$ to some other vertex $y_2$.
Then they form a $\Theta_{p^2-1}$ inside $G$, which is also impossible.
In the next subcase, the edges begin at one common vertex $y_0$
and end at $p^2$ distinct vertices $y_1, \ldots, y_{p^2}.$  Since the
graph $G$ is admissible, it has no free edges and every nonbasepointed
vertex has valence at least $3$.  So all of the vertices
$y_1, \ldots, y_{p^2}$ have to connect up in some manner, and in doing
so they will violate the fact that $\pi_1(G)=F_{2p-1}$.  The final
three subcases, in which the $e_i$ either form a $p^2$-gon, have no
common vertices,
or form loops with $p^2$ distinct endpoints,
are similar.  Hence $P$ will never
be $\Z/p^2$ and this case will not occur.

\medskip

\noindent CASE 2.  Next we will examine the case where
$P=\Z/p \times \Z/p= (\alpha) \times (\beta)$. 
The first cyclic summand must rotate $p$ edges
$e_1, e_2, \ldots, e_p$ of $\Gamma$. 
Without loss of
generality, we may assume that the basepoint $*$ is one of the
endpoints of each $e_i$.  Now if $\beta$ sends all of the $e_i$ to
another whole collection 
$\beta e_i$ (with $\{ e_i \}$ disjoint from $\{ \beta e_j \}$) then
the basepoint $*$ must be one of the endpoints of each $\beta^j e_i$ also; 
therefore, we obtain at least $p^2$ edges emanating from the basepoint $*$
which are moved by $\alpha$ and $\beta$.  This implies that the rank of
$\pi_1(\Gamma)$ is at least $p(p-1)$ (i.e., the best that can happen is
that $p$ copies of $\Theta_{p-1}$ are wedged together at the basepoint),
which is too large as $p \geq 5$.

So $\beta$ does not send the $e_i$ to another whole collection
$\beta e_i$ of edges disjoint from the $e_i$.
Without loss of generality,
$(\beta)$ fixes the edges $e_i$ (by replacing $\beta$ by
$\beta - \alpha^j$ if necessary.)  Hence the collection $\{ e_i \}$
is $P$-invariant. Now $\beta$
must rotate $p$ other edges $f_1, f_2, \ldots, f_p$.  
The $e_i$ do not form a $p$-gon as $*$ is an endpoint of each of 
them.  Hence the $e_i$ form either a rose, a star, 
or a $\Theta$-graph. If they form a rose or a star,
then the $f_i$ must form a $\Theta_{p-1}$ else the
rank of $\pi_1(\Gamma)$ is larger than $2p-1$.  If the $e_i$ form a 
$\Theta$-graph $\Theta_{p-1}$, then there are $p$ holes available
in the rank of $\pi_1(\Gamma)$ for
the other edges of $\Gamma$ to use up. 

By doing the sort of case-by-case analysis
that we did in the previous paragraph, we see that $G$ must be one of the
following graphs (listed in increasing order with respect
to the number of vertices):
\begin{itemize}
\item $R_p \vee \Theta_{p-1}$, whose automorphism group has
the same cohomology as $\Sigma_p \times \Sigma_p.$
\item $\Theta_{2p-1}$, whose
automorphism group is $\Sigma_{2p}$.
\item $\Theta_{p-1} \vee \Theta_{p-1}$, plus one additional edge $e$
attached in some manner to the existing vertices.
The automorphism group here will have the same
cohomology as either
$\Sigma_p \times \Sigma_p$ or 
$(\Sigma_p \times \Sigma_p) \rtimes \Z/2$.
\item $\Theta_{p-1} \coprod \Theta_{p-1}$,
with one additional
edge $e_1$ attached going from an already existing vertex
of one of the $\Theta$-graphs to one on the other $\Theta$-graph,
after which we sequentially attach another
edge $e_2$ to that resulting graph.
The endpoints of $e_2$ can be attached to any of the
already existing vertices, or they can be attached anywhere in
the interior of $e_1$. 
The automorphism group here will have the same
cohomology as either
$\Sigma_p \times \Sigma_p$ or 
$(\Sigma_p \times \Sigma_p) \rtimes \Z/2$.
\item $\Theta_{p-1} \vee \Xi_p$, with
automorphism group $\Sigma_p \times D_{2p}$.
\end{itemize}

\medskip

\noindent CASE 3.  For the final case, $P=\Z/p=(\alpha)$.
We want to show that all of the $p$-symmetries in
the graph $G$ are also at least $D_{2p}$-symmetries.  That is, in addition
to the rotation by $\Z/p$, there is also a
dihedral ``flip''.
We will be
able to get this result because $n=2p-1$ is not large enough with
respect to $p$ for us to be able to generate graphs $G$ with
$\pi_1(G)=F_n$ that have $\Z/p$-symmetries but not
$D_{2p}$-symmetries.

We have $P = \Z/p$
acting on a graph whose fundamental group has
rank $n=2p-1$.
As before, there exist at least $p$
edges $e_1, \ldots, e_p$ that $P$ rotates.  If these edges form a
$\Theta_{p-1}$ or an $R_p$, we are done.  This is because now
$P$ cannot move any other edges of $G$, else we are in the case of the
previous paragraph where  $P=\Z/p \times \Z/p$.
As the automorphism groups of both $\Theta_{p-1}$
and $R_p$ contain the symmetric group $\Sigma_p$,
we are done in this subcase.

Now suppose we are in the other extreme subcase, the one where
$e_1, \ldots, e_p$ have no endpoints in common.  Choose a minimal path
$\gamma_1$ from $e_1$ to the basepoint.  Since $P=(\alpha)$ fixes the
basepoint, we have that $\gamma_i := \alpha^{i-1} \gamma_1$ is a
minimal path from $e_i$ to $x_0$ for all $i=1, \ldots, p$.  Since none
of the endpoints of the $e_i$ can be the basepoint, there are at least
$p$ distinct edges $f_1, \ldots, f_p$ in $\gamma_1, \ldots, \gamma_p$,
respectively.  We can also assume that each $f_i$ has at least one
endpoint that is not the basepoint. Because $G$ has no free edges,
no separating edges, 
and all non-basepointed vertices have valence at least three, another
case-by-case analysis reveals that
since $\pi_1(G)$ must have rank less than $2p$, the graph is forced
to be either $\Upsilon_{2p-1}^1$, $\Upsilon_{2p-1}^2$,
or $\Lambda_{2p-1}^1$.  The first two of these graphs have
dihedral symmetry, while the last has
automorphism group with the same cohomology as $\Sigma_p$.

The next case is the one in
which the $e_i$ are all loops with $p$ distinct
endpoints $y_i$.  As in the previous case, we can choose
a $\Z/p$ equivariant path from each $y_i$ to the
basepoint.  The admissibility conditions on the graph
only allow one possibility, namely the graph $\Lambda_{2p-1}^2$.
As this graph has automorphism group with the same
cohomology as $\Sigma_p$, we are finished with this case.

Next consider the case where the $e_i$ form a $p$-gon.  Since the
vertices of $G$ have valence at least $3$, there must be $p$ other
edges $f_1, \ldots, f_p$ in $G$ that each start at one of the $p$
vertices of the $p$-gon.  Since $G$ is admissible and the rank of
$\pi_1(G)$ is $2p-1$, these additional edges cannot also join up to
form a $p$-gon.  (Why?  Both $p$-gons still need to connect up to the
basepoint in some way, and in connecting up to $x_0$ the rank of the
fundamental group of $G$ will be forced too high.)  In addition, the
edges must have some vertices in common, else we reduce to the
previous case; therefore, the $f_i$ are all forced to end at some common
vertex $y_0$.  In other words, we have a $\Xi_p$ embedded in $G$.  If
$P$ doesn't move any other edges in $G$, we are done since $\Xi_p$
has dihedral symmetry.  If some other edges $g_1, \ldots, g_p$
are moved, they must also
be attached to the $p$-gon that the $e_i$ form, or we will have two
independent $\Z/p$-actions and be in the case $P=\Z/p \times \Z/p$.
None of the following cases can happen, else
rank$(\pi_1(G)) \geq 2p$:
\begin{itemize}
\item The other endpoints of the $g_i$ all connect to $y_0$.
\item The other endpoints of the $g_i$ also connect to the $p$-gon
formed by the $e_i$.
\item The other endpoints of the $g_i$ form $p$ other distinct
vertices.
\end{itemize}
Hence these other endpoints all have to connect to some other common
vertex $y_1$, forming another copy of $\Xi_p$ in $G$.  Thus $G$ must
be the graph $\Upsilon_{2p-1}$, which certainly has dihedral symmetry.

\input{fil1p.pic}
\begin{figure}[here]
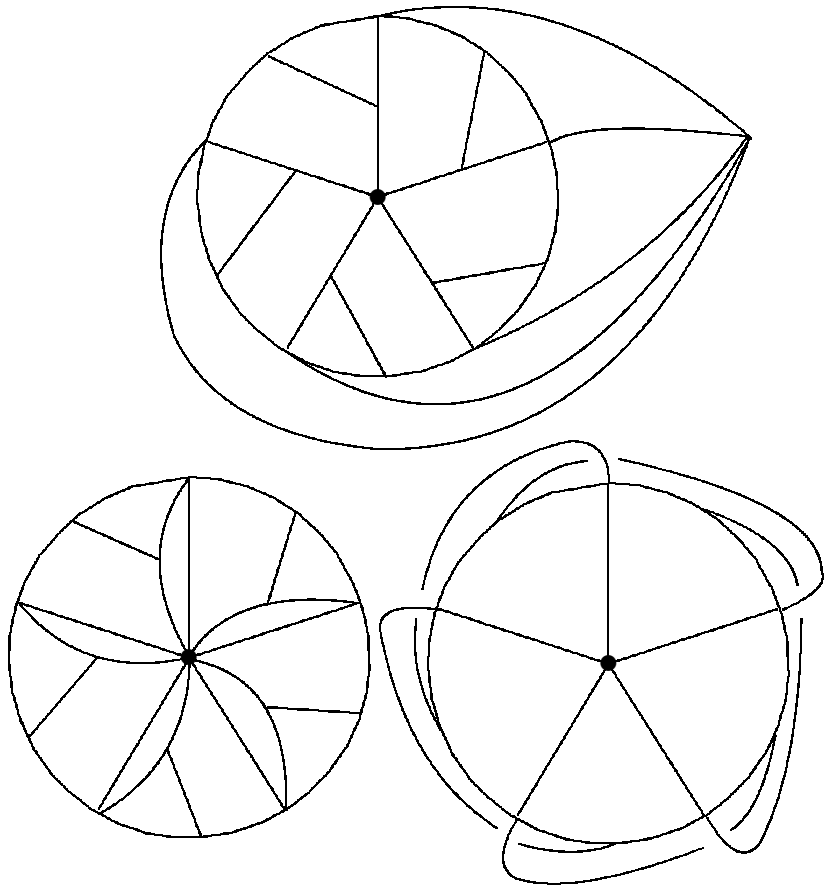

\caption{\label{fil1} Graphs whose symmetry groups are exactly $\Z/p$}
\end{figure}

For the final case, the edges $e_i$ have one common vertex $y_0$,
and end in $p$ other distinct vertices $y_1, \ldots, y_p$.
In addition, i) there are no $p$-gons in $G$, ii) there
are no collections
of $p$ edges in $G$ that are rotated by $P$ and that have no common
vertices,
and iii) there
are no collections
of $p$ edges in $G$ that are rotated by $P$ and that
each form loops with distinct endpoints.
Since all of $y_1, \ldots, y_p$ have valence three, $P$
must rotate two other collections of edges
$\{f_1, \ldots, f_p\}$ and $\{g_1, \ldots, g_p\}$ that begin at the
vertices $y_1, \ldots, y_p$ and end at the vertices $z_1$ and
$z_2$, respectively.  Note that since $\pi_1(G)=F_{2p-1}$,
$|\{y_0,z_1,z_2\}| \geq 2$.  Also note that $P$ cannot move any other
edges of $G$ except the ones we have listed.  In this case, the
symmetric group $\Sigma_p$ acts on the collections of edges defined above,
and so the cohomology of the group of graph automorphisms of
the graph is the same as that of the symmetric group.
If $|\{y_0,z_1,z_2\}| = 2$ then the only edges in the
graph are the $e_i$, $f_i$, and $g_i$ and the graph is
either $\Upsilon_{2p-1}^{2a}$ or $\Upsilon_{2p-1}^{2b}$.
On the other hand, if $|\{y_0,z_1,z_2\}| = 3$, then
the graph has one additional edge besides the
$e_i$, $f_i$, or $g_i$.
Accordingly, the graph looks like
a $\Phi_{2p-1}$ (See Figure \ref{fig1}) with
one additional edge added.  This additional edge can go
from any of the $\{y_0,z_1,z_2\}$ to any other one,
including possibly the
same one.  In any case, it is definitely true that the
graph has automorphism group with the same
cohomology as $\Sigma_p$.
The lemma follows. \end{proof}

For an example of what we were trying to avoid in the proof
of the above lemma,
refer to the three examples given in the Figure \ref{fil1}.
The graphs pictured have an obvious $\Z/p$-symmetry
given by rotation about the basepoint, which is
indicated by a solid  dot.  But they have no dihedral
flip, and their basepoint-preserving automorphism groups
are all exactly $\Z/p$,
where $p=5$ in the examples pictured and where obvious
analogues exist for other odd primes.  The ranks of the
fundamental groups of the graphs pictured are
$3p-1$, $3p$, and $2p$, respectively.  The last rank, $2p$,
is the lowest rank possible where one can have a graph with
exactly $\Z/p$ symmetry. 


\begin{cor} \label{ttt2}
A vertex in the $p$-singular locus of $X_n$ has at most dihedral symmetry 
if its cohomology is the same as that of $D_{2p}$ or $D_{2p} \times \Sigma_p$,
and vertices in the $p$-singular locus
will never have exactly $\Z/p$ symmetry.
\end{cor}

In the figures below, a dotted line or a hollow dot indicates
that the given edge or vertex, respectively, does not have
the indicated property.  A solid dot, a solid line, or
a $2$-simplex with an X in it, means that the given
vertex, edge, or $2$-simplex, respectively, does have
the indicated property.

By analyzing the $\Z/p$-invariant subforests of all of the graphs
explicitly listed in the proof of Lemma \ref{t40}, we can see what
types of stabilizers higher dimensional simplices
(rather than just vertices) have.

\begin{figure}[here]
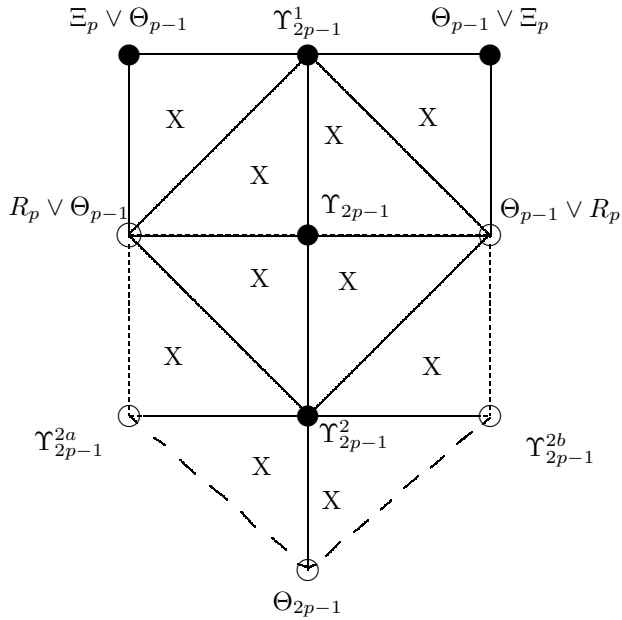

\caption{\label{fig4} Some simplices with at most dihedral symmetry}
\end{figure}

We will show that the
simplices with at most dihedral symmetry will fall into
two (exhaustive but not disjoint) categories.  The first
category
consists of those that are listed in
Figure \ref{fig4}.
The second category consists of simplices whose maximal vertex
(recall that $X_n$ is the realization of a poset)
has the form $\Xi_p \vee \Gamma_{p-1}$ 
where $\Gamma_{p-1}$ is some basepointed
graph with fundamental group of rank $p-1$, the wedge does
not necessarily take place at the basepoint, and where the
forest collapses of the simplex respect the $\Z/p$ action
on $\Xi_p$.

\bigskip

\input{fig4p.pic}

\bigskip

\input{fil2p.pic}
\begin{figure}[here]
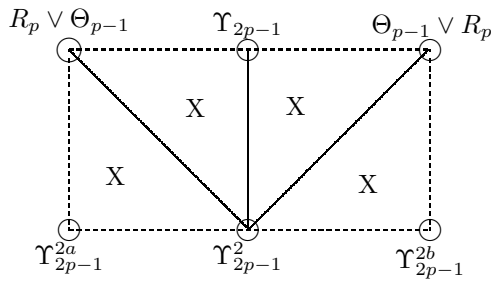

\caption{\label{fil2} Simplices with exactly $\Z/p$ symmetry}
\end{figure}

We will also show that the simplices listed in
Figure \ref{fil2} are the only ones with exactly $\Z/p$
symmetry.

\begin{cor} \label{t41} Let $p \geq 5$ be prime, $n = 2p-1$,
and consider the $p$-singular locus of the spine $X_n$ of
auter space.
\begin{itemize}
\item
The only simplices
with at most dihedral symmetry are either: (i) listed in
Figure \ref{fig4};  or (ii) have maximal vertex
of the form $\Xi_p \vee \Gamma_{p-1}$.
\item The only simplices with exactly
$\Z/p$ symmetry are those listed in Figure \ref{fil2}.
\end{itemize}
\end{cor}

\begin{proof} 
We examine each of the graphs listed in
Lemma \ref{t40} separately.  By enumerating the $\Z/p$ invariant
subforests of each of these graphs, one can list all of the simplices
in the $p$-singular locus of $X_n$.  We can ignore the graphs
in Lemma \ref{t40} that do not have dihedral symmetry, as all of
their symmetry comes from symmetric groups.  When you collapse
invariant subforests of these graphs, you still get
graphs with symmetry coming from the symmetric group.

So we are left with analyzing the graphs from Lemma \ref{t40} with
dihedral symmetry, which were:
\begin{itemize} 
\item $\Theta_{p-1} \vee \Xi_p$. (There are actually two
possibilities here as the enumeration in Lemma \ref{t40}
did not specify basepoints.  The central vertex of $\Xi_p$
could be attached to either the basepoint of $\Theta_{p-1}$
or the other vertex of $\Theta_{p-1}$.)
\item $\Gamma_{p-1} \vee \Xi_p$, where $\Gamma_{p-1}$
is a basepointed graph with fundamental group of rank $p-1$
which has
no $p$-symmetry (or where $\Gamma_{p-1}$ is $\Theta_{p-1}$
but the central vertex of $\Xi_p$ is attached to the
midpoint of an edge of $\Gamma_{p-1}$.)
\item $\Upsilon_{2p-1}$.
\item $\Upsilon_{2p-1}^1$.
\item $\Upsilon_{2p-1}^2$.
\end{itemize}

For the first two types of graphs, you can obtain
simplices with dihedral symmetry by collapsing all of the spokes
of $\Xi_p$ and/or any forest in the other graph
of the wedge sum
(either $\Theta_{p-1}$ or $\Gamma_{p-1}$.)  The resulting
simplex with maximal vertex $\Theta_{p-1} \vee \Xi_p$
or $\Gamma_{p-1} \vee \Xi_p$ will
clearly have at most dihedral symmetry, and 
will also just as clearly not give you a graph with exactly $\Z/p$
symmetry.

In a similar manner, simplices in the $p$-singular locus
of $X_n$ with maximal vertex $\Upsilon_{2p-1}$ or
$\Upsilon_{2p-1}^1$ are (exhaustively) listed in
Figure \ref{fig4}.  Note that $\Upsilon_{2p-1}$
can only be blown up (while still preserving the $\Z/p$
action so that we stay in the $p$-singular locus of $X_n$)
in two ways, to either
$\Upsilon_{2p-1}^1$ or $\Upsilon_{2p-1}^2$. The latter two
graphs cannot be blown up at all.

Finally, the simplices with maximal vertex  $\Upsilon_{2p-1}^2$
are listed in Figure \ref{fig4} or Figure \ref{fil2}.
Note that we can obtain edges and $2$-simplices with
exactly $\Z/p$ symmetry, even though no actual vertex
of $X_n$ has exactly $\Z/p$ symmetry.  This is because
you can choose subforests of $\Upsilon_{2p-1}^2$
that do respect the dihedral ``flip'' of
$\Upsilon_{2p-1}^2$.  In other words, this flip will
not take the subforest to itself again.  Hence the
resulting simplex will just have symmetry group $\Z/p$.
Last of all, note that you can also choose subforests
of $\Upsilon_{2p-1}^2$ which do respect the dihedral
flip, and these give simplices with dihedral symmetry. \end{proof}

\section{The integral cohomology of the quotient never stabilizes} \label{c14}

We will prove Theorem \ref{t44} in this section.
As in Section \ref{c13}, all primes $p$
considered are assumed to be greater than or equal to $5$.

\begin{lemma} \label{t42}
For the rows $0 \leq s < 2(p-1)$, the $E_2$ page of
the spectral sequence \ref{e1}
applied to calculate $H^*(Aut(F_{2p-1}); \Z_{(p)})$
is given by
\bigskip
$$\matrix{
\hfill E_2^{r,s} &=& &\left\{\matrix{
H^r(Q_{2p-1}; \Z_{(p)})  \hfill &s = 0 \hfill \cr
\Z/p \hfill &r=2 \hbox{ and } s = 4k-2 > 0, k \in \Z^+ \hfill \cr
0 \hfill &\hbox{otherwise} \hfill \cr} \right. \hfill \cr
}$$
\end{lemma}
\bigskip

\begin{proof} As $E_1^{r,0}$ is the cochain complex $C^r(Q_{2p-1}; \Z_{(p)})$,
it follows that $E_2^{r,0} = H^r(Q_{2p-1}; \Z_{(p)})$ as
claimed above.

None of the simplices in $X_{2p-1}$ contribute anything to the
odd rows between $0$ and $2(p-1)$ of the above spectral sequence,
from Corollary \ref{t41}.  Also from Corollary \ref{t41}, the ones that contribute
to rows of the form $4k - 2$,
$k \in \Z^+$, are all listed in Figure \ref{fil2}.
Let $A$ be the subcomplex of $Q_n$ generated by all of the
simplices pictured in Figure \ref{fil2} and let $B$ be the
subcomplex generated by just the simplices corresponding to
dotted lines or hollow dots in Figure \ref{fil2}.  Then the row
$s=4k-2$ on the $E_1$ page of the spectral sequence is
$C^r(A,B; \Z/p)$.  Examining Figure \ref{fil2} we see that
$$\matrix{
\hfill H^r(A,B; \Z/p) &=& &\left\{\matrix{
\Z/p \hfill &r=2 \hfill \cr
0 \hfill &\hbox{otherwise} \hfill \cr} \right. \hfill \cr
}$$
Consequently the $E_2$ page is as claimed for the
rows $s = 4k-2$.

Our final task is to calculate the $E_2$ page for the
rows $s=4k$.  Simplices in the $p$-singular locus
of the spine with ``at most dihedral symmetry'' contribute
to these rows.  From Corollary \ref{t41}, we have a characterization
of such simplices. 
Define the subcomplex $M$ of the
$p$-singular locus of the spine
$X_{2p-1}$ of auter space to be the
subcomplex generated by
simplices with ``at most dihedral symmetry''.
More precisely, from Corollary \ref{t41}, we know it is
generated by the simplices corresponding
to those in Figure \ref{fig4}
(i.e., corresponding in the sense that we are
taking $M$ to be a subcomplex of the spine
rather than its quotient and Figure \ref{fig4} is
a picture in the quotient)
in addition to simplices whose maximal vertex
has underlying graph of the
form $\Xi_p \vee \Gamma_{p-1}$
(where the
forest collapses in the simplices respect
the $\Z/p$ action on $\Xi_p$.)
Recall that an $r$-simplex 
with at most dihedral symmetry
contributes exactly one $\Z/p$ to 
$E_1^{r,4k}$, while all other
simplices (those without dihedral or
exactly $\Z/p$ symmetry)
contribute nothing to this row.
 
Let $N$ be the
subcomplex of $M$ generated by simplices in $M$
which do not have at most dihedral symmetry.
Observe that none of the simplices in $N$ have
at most dihedral symmetry.
Also note that the row $E_1^{*,4k}$ is the
relative cochain complex
$$C^*(M/Aut(F_{2p-1}),N/Aut(F_{2p-1}); \Z/p).$$

Let $M'$ be the subcomplex of $M$ generated by $N$ and by
simplices whose maximal vertex is 
$\Upsilon_{2p-1}^2$.  Hence $M'$ is the subcomplex
consisting of $N$ and the bottom two thirds
of Figure $\ref{fig4}$.  There is an 
$Aut(F_{2p-1})$-equivariant deformation retraction
of $M$ onto $M'$, given on the vertices of the
poset by:
\begin{itemize}
\item Contracting the spokes of the
graph $\Xi_p$ in $\Theta_{p-1} \vee \Xi_p$.
\item Contracting the spokes of the graph $\Xi_p$ in
$\Gamma_{p-1} \vee \Xi_p$, where $\Gamma_{p-1}$ has
no $p$-symmetry.
\item Contracting the $p$ outward radiating
edges attached to the $p$-gon in the center of
the graph $\Upsilon_{2p-1}^1$.  In the terminology
used at the beginning of Section \ref{c13} while defining
$\Upsilon_{2p-1}^1$, we are contracting the edges $e_i$.
\end{itemize}
That it is a deformation retraction follows from the
Poset Lemma in \cite{[K-V]} attributed to Quillen.

%

As the homotopy retracting $M$ to $M'$ is
$Aut(F_{2p-1})$-invariant, it descends to a deformation
retraction of $M/Aut(F_{2p-1})$
to $M'/Aut(F_{2p-1})$.  Hence the relative cohomology
groups
$$H^*(M/Aut(F_{2p-1}),N/Aut(F_{2p-1}); \Z/p)$$
and
$$H^*(M'/Aut(F_{2p-1}),N/Aut(F_{2p-1}); \Z/p)$$
are isomorphic.  Now referring to
Figure $\ref{fig4}$,
we see that
$$H^t(M'/Aut(F_{2p-1}),N/Aut(F_{2p-1}); \Z/p) = 0$$
for all $t$
because we can contract all of the
simplices in $M'/Aut(F_{2p-1})$
uniformly into $N/Aut(F_{2p-1})$. \end{proof}

An immediate consequence is


\begin{proof}[Proof of Theorem \ref{t44}:]
From \cite{[H-V]}, if $m \geq 8k+3,$ then the standard map
$$H^{4k}(Aut(F_{m+1}); \Z) \to H^{4k}(Aut(F_m); \Z)$$
is an isomorphism.
Observe that $H^{4k}(Aut(F_{8k+3}); \Z) =
H^{4k}(Aut(F_\infty); \Z)$ is a finitely generated
abelian group.  If it contains a torsion free summand isomorphic to
$\Z$, then we are done and $H^{4k}(Aut(F_\infty); \Q) \not =0$.
Otherwise, choose a prime $q$ such that $2q-1 \geq 8k+3$
and so that for all primes $p \geq q$ there is no $p$-torsion
in $H^{4k}(Aut(F_{8k+3}); \Z)$.
We will show that $H^{4k+1}(Q_{2p-1}; \Z)$  has $p$-torsion for all
primes $p \geq q$, which will prove the theorem.

Let $p \geq q$.  From the lemma above, if we
use the standard equivariant spectral sequence to calculate
$H^{*}(Aut(F_{2p-1}); \Z_{(p)})$, then a class
$\alpha \in E_1^{2,4k-2}$ in the $E_1$-page survives at least
until the $E_{4k-1}$-page.

Because $H^{4k}(Aut(F_{2p-1}); \Z)$ has no
$p$-torsion and $H^{4k}(Aut(F_{2p-1}); \Q) = 0$, 
we have 
$H^{4k}(Aut(F_{2p-1}); \Z_{(p)})=0$.  
Hence the class $\alpha \in E_1^{2,4k-2}$ cannot survive to the
$E_\infty$ page. 
It follows that there is $p$-torsion in
$E_{4k-1}^{4k+1,0}$. 
Recall that
$E_1^{r,0}$ corresponds to the cellular
chain complex with $\Z_{(p)}$ coefficients
for $Q_{2p-1}$.  The $p$-torsion in
$E_{4k-1}^{4k+1,0}$, therefore, would have to have been created
when going from the $E_1$ to $E_2$ pages, because any of
the torsion above the horizontal axis of the spectral sequence
could not map onto a torsion free element on the horizontal axis.
So $H^{4k+1}(Q_{2p-1}; \Z_{(p)})$ has $p$-torsion, and thus
$H^{4k+1}(Q_{2p-1}; \Z)$ has $p$-torsion. \end{proof}


\begin{thebibliography}{1}

\bibitem{[B-T]}
G, Baumslag and T. Taylor,
{\em The center of groups with one defining relator},
Math. Ann. 175 (1968) 315-319.

\bibitem{[B]}
K. Brown,
{\em Cohomology of Groups},
Springer-Verlag Berlin, Heidelberg 1982.

\bibitem{[C-V]}
M. Culler and K. Vogtmann,
{\em Moduli of graphs and automorphisms of free groups},
Invent. Math. 84 (1986) 91-119.

\bibitem{[G-M]}
H. H. Glover and G. Mislin,
{\em On the $p$-primary cohomology of $Out(F_n)$ in the $p$-
rank one case},
to appear in J. Pure Appl. Algebra 150 (2000).

\bibitem{[H]}
A. Hatcher,
{\em Homological stability for automorphism groups of free groups},
Comment. Math. Helv. 70 (1995) 39-62.

\bibitem{[H-V]}
A. Hatcher and K. Vogtmann,
{\em Cerf theory for graphs},
J. London Math. Soc. (2) 58 (1998) 633-655.

\bibitem{[V]}
A. Hatcher and K. Vogtmann,
{\em Rational homology of $Aut(F_n)$},
Math. Res. Lett. 5 (1998), no. 6, 759-780.

\bibitem{[JD]}
C. A. Jensen,
{\em Cohomology of $Aut(F_n)$},
Cornell University Ph.D. dissertation, Ithaca, New York 1998.

\bibitem{[J1]}
C. A. Jensen,
{\em Cohomology of $Aut(F_n)$ in the $p$-rank two case},
to appear in J. Pure Appl. Algebra.

\bibitem{[K-V]}
S. Krstic and K. Vogtmann,
{\em Equivariant outer space and automorphisms of
free-by-finite groups},
Comment. Math. Helv. 68 (1993) 216-262.

\bibitem{[M]}
H. Minkowski,
{\em Zur Theorie der positiven quadratischen Formen},
Crelles J. 101 (1887) 196-202.

\bibitem{[S-V1]}
J. Smillie and K. Vogtmann,
{\em A generating function for the Euler characteristic
of $Out(F_n)$},
J. Pure Appl. Algebra 44 (1987) 329-348.

\bibitem{[S-V2]}
J. Smillie and K. Vogtmann,
{\em Automorphisms of graphs, $p$-subgroups of $Out(F_n)$ and the
euler characteristic of $Out(F_n)$},
J. Pure Appl. Algebra 49 (1987) 187-200.

\bibitem{[S]}
R. G. Swan,
{\em The $p$-period of a finite group},
Ill. J. Math 4 (1960), 341-346.

\end{thebibliography}
\end{document}